\documentclass[a4paper,11pt]{article}

\usepackage{amsmath,amssymb,amsthm,ifthen,lscape,setspace}
\usepackage{graphicx}
\usepackage[curve,knot]{xypic}

\usepackage[a4paper,left=1in,width=458pt,textheight=650pt,top=1.25in]{geometry}

\usepackage{stmaryrd,bbm}

\usepackage{hyperref}

\theoremstyle{plain}
\newtheorem{theorem}{Theorem}[section]
\newtheorem*{theorem*}{Theorem}
\newtheorem{proposition}[theorem]{Proposition}
\newtheorem{lemma}[theorem]{Lemma}
\newtheorem{corollary}[theorem]{Corollary}
\theoremstyle{definition}
\newtheorem{definition}[theorem]{Definition}
\theoremstyle{remark}

\newtheorem*{example*}{Example}
\newtheorem*{examplectd*}{Example (continued)}

\newenvironment{jegrmk}{\vspace{0.5em}\par \noindent \textit{Remark.}}{\vspace{0.5em}}

\newcounter{property}
\newcounter{subproperty}
\newenvironment{property}{\refstepcounter{property} \par \vspace{1em} \noindent (\arabic{property})\begin{math}\quad }{\end{math} \vspace{1em} \\}

\input cyracc.def

\newcommand{\keywords}[1]{\noindent \footnotesize {\textbf{Keywords:} #1 } \normalsize }

\newboolean{Sectext}
\newsavebox{\MSCtextP}
\newsavebox{\MSCtextS}
\newcommand{\MSC}[2]{
\sbox{\MSCtextP}{#1}
\ifthenelse{\equal{#2}{NONE}}{\setboolean{Sectext}{false}}{\setboolean{Sectext}{true}}
\sbox{\MSCtextS}{#2}
}
\newcommand{\printMSCwithtitle}{
\ifthenelse{\boolean{Sectext}}{\noindent \footnotesize \textbf{Mathematics Subject Classification (2000):} \usebox{\MSCtextP} (Primary), \usebox{\MSCtextS} \normalsize}{\noindent \footnotesize \textbf{Mathematics Subject Classification (2000):} \usebox{\MSCtextP} \normalsize}
}

\newcommand{\ad}[2]{{\mathrm{ad}}_{#1}(#2)}

\newcommand{\bigdsum}{\bigoplus}

\newcommand{\bdelta}{\underline{\;}\mkern-6mu{\delta}}

\let\chisave\chi
\renewcommand{\chi}{{%
 \mathchoice{\raisebox{0.25ex}{$\displaystyle\chisave$}}
            {\raisebox{0.2ex}{$\textstyle\chisave$}}
            {\raisebox{0.2ex}{$\scriptstyle\chisave$}}
            {\raisebox{0.1ex}{$\scriptscriptstyle\chisave$}}}}

\newcommand{\cross}{\times}

\newcommand{\dbos}[3]{\Lie{#1} \rtimesdot \Lie{#2} \ltimesdot \Lie{#3}}

\newcommand{\defeq}{\stackrel{\scriptscriptstyle{\mathrm{def}}}{=}}
\newcommand{\dform}{\mathrm{d}}
\newcommand{\dsum}{\ensuremath{ \oplus}}
\newcommand{\dual}[1]{\ensuremath {#1}^{*}}

\renewcommand{\epsilon}{\varepsilon}

\newcommand{\id}{\ensuremath \mbox{\textup{id}}}

\newcommand{\inj}{\hookrightarrow}
\newcommand{\integ}{\ensuremath{\mathbb{Z}}}
\newcommand{\intersection}{\mathrel{\cap}}
\newcommand{\inv}[1]{\ensuremath {#1}^{-1}}
\newcommand{\ip}[2]{\ensuremath \lgen\;\!#1,#2\;\!\rgen}

\newcommand{\iso}{\ensuremath \cong}

\newcommand{\Ker}[1]{\ensuremath \mbox{Ker}\>{#1}}
\newcommand{\laction}{\triangleright}

\newcommand{\lgen}{\ensuremath \mathopen{<}} 

\newcommand{\ltimesdot}{\mathrel{\cdot\joinrel\mkern-14.3mu\rhd\mkern-8.7mu<}}
\newcommand{\Lbracket}[2]{\ensuremath{ [\, #1 , #2\, ]}}
\newcommand{\Lbbracket}[2]{\ensuremath{ \llbracket\, #1 , #2\, \rrbracket}}
\newcommand{\Lie}[1]{\ensuremath{\mathfrak{#1}}}

\newcommand{\modcat}[1]{\ensuremath \!_{#1}\mathcal{M}}

\newcommand{\nat}{\ensuremath \mathbb{N}}

\newcommand{\onto}{\twoheadrightarrow}
\newcommand{\op}[1]{{#1}^{\mbox{\scriptsize \textup{op}}}}

\renewcommand{\phi}{\varphi}

\newcommand{\rank}{\mathop{\mathrm{rank}}}

\newcommand{\rgen}{\ensuremath \mathclose{>}}

\newcommand{\rtimesdot}{\mathrel{>\joinrel\mkern-6mu\lhd\mkern-7.2mu\cdot}}

\newcommand{\smun}[2]{{#1}_{\scriptscriptstyle #2}}

\newcommand{\tensor}{\ensuremath \otimes}

\newcommand{\union}{\mathrel{\cup}}
\newcommand{\ve}[1]{{#1}^{\scriptscriptstyle{\vee}}}

\newcommand{\cf}{cf.\ }
\newcommand{\ie}{i.e.\ }
\newcommand{\eg}{e.g.\ }
\newcommand{\webtilde}{\footnotesize{$\sim\!$} \small}

\title{Braided-Lie bialgebras associated to Kac--Moody algebras}
\author{Jan E. Grabowski\footnotemark[2] 
\\ \small{\textit{Keble College, Oxford, OX1 3PG, United Kingdom}}}
\date{\today}

\hypersetup{
pdftitle={Braided-Lie bialgebras associated to Kac-Moody algebras},
pdfauthor={Jan E. Grabowski},
pdfstartview=FitH,
pdfkeywords={Kac-Moody algebra, braided, Lie bialgebra}
}

\begin{document}

\maketitle

\renewcommand{\thefootnote}{\fnsymbol{footnote}}
\footnotetext[2]{Email: \href{mailto:jan.grabowski@maths.ox.ac.uk}{jan.grabowski@maths.ox.ac.uk}.  Website: \href{http://www.maths.ox.ac.uk/~grabowsk/}{http://www.maths.ox.ac.uk/\webtilde \footnotesize grabowsk/ \normalsize}}
\renewcommand{\thefootnote}{\arabic{footnote}}
\setcounter{footnote}{0}

\vspace{-1.15em}
\begin{abstract}
\noindent Braided-Lie bialgebras have been introduced by Majid, as the Lie versions of Hopf algebras in braided categories.  In this paper we extend previous work of Majid and of ours to show that there is a braided-Lie bialgebra associated to each inclusion of Kac--Moody bialgebras.  Doing so, we obtain many new examples of infinite-dimensional braided-Lie bialgebras.  We analyze further the case of untwisted affine Kac--Moody bialgebras associated to finite-dimensional simple Lie algebras.  The inclusion we study is that of the finite-type algebra in the affine algebra.  This braided-Lie bialgebra is isomorphic to the current algebra over the simple Lie algebra, now equipped with a braided cobracket.  We give explicit expressions for this braided cobracket for the simple Lie algebra $\mathfrak{sl}_{3}$.
\end{abstract}

 \\
\MSC{\footnotesize 17B67}{\footnotesize 17B62, 22E67}
\printMSCwithtitle

\section{Introduction}\label{s:intro}
Following Drinfel\cprime d's ICM address (\cite{DrinfeldICM}) some twenty years ago, there has been much interest among mathematicians and physicists in the extension of the classical theory of Lie algebras by considering not just a Lie bracket but in addition a compatible Lie cobracket.  Objects with these structures are called Lie bialgebras and play the same role as Hopf algebras do with respect to algebras.  Many people have studied Lie bialgebra structures on the well-known classes of Lie algebras, not least because in many cases these yield solutions to the classical Yang--Baxter equation and hence have connections to integrable systems and knot theory.

Another direction in Hopf algebra theory has been emphasized by Majid, namely the study of braided monoidal categories and in particular algebraic structures on objects in such categories.  By ``in'' here we mean not just membership of the category but also covariance of the algebraic structure in question with respect to the braiding.  Typically the braided monoidal category is a module category of some sort and then we mean covariance with respect to the action.  The books \cite{FQGT} and \cite{QGP} describe much of the theory of algebraic structures in braided categories developed up to the start of this millennium.  More recent work has included that of Andruskiewitsch and Schneider (\cite{AndrSchnPHA}) and others on pointed Hopf algebras and Nichols algebras, the latter being certain types of Hopf algebras in braided categories.

Majid (\cite{BraidedLie}) has also explained how the ideas of braided categories apply in the Lie algebra case and has introduced the Lie analogue of Hopf algebras in braided categories, which he called braided-Lie bialgebras.  These objects are the focus of our work here: we show how they arise very naturally and extend Majid's treatment to produce many new infinite-dimensional examples.  We do this by extending previous results of Sommerh\"{a}user (\cite{SommerhauserKacMoody}), Majid (op.\ cit.) and ourselves (\cite{LieInduction}) to the case of Kac--Moody algebras, the class of Lie algebras associated to generalized Cartan matrices discovered independently by Kac (\cite{Kac1967}) and Moody (\cite{Moody1967}) in the mid-1960s.  These Lie algebras are typically infinite-dimensional but retain many of the nice structural features of the finite-dimensional semisimple Lie algebras.

The first aim of this work is to show that braided-Lie bialgebra structures are naturally associated to inclusions of Kac--Moody algebras.  The second is to provide a more detailed analysis of this situation when the larger algebra is the untwisted affine Kac--Moody algebra over a finite-type simple algebra, embedded in degree $0$.  In this case, we show that the associated braided-Lie bialgebra is isomorphic to the current algebra over the simple subalgebra, now equipped with a braided cobracket. 

Also of interest is Majid's double-bosonisation construction (\cite{BraidedLie}), which takes as input a quasitriangular Lie bialgebra $\smun{\Lie{g}}{\mathrm{0}}$ and a braided-Lie bialgebra $\Lie{b}$.  One then obtains a new quasitriangular Lie bialgebra $\Lie{g} \iso \dbos{b}{\smun{g}{\mathrm{0}}}{\op{\dual{b}}}$.  (Strictly, this is only true if $\Lie{b}$ is finite-dimensional: in general, one should take two dually paired braided-Lie bialgebras.)  Here $\rtimesdot$ is a simultaneous semidirect Lie algebra and coalgebra structure; similarly for $\ltimesdot$ on the other side.  It is possible to realize this as a quotient of the Drinfel\cprime d double (\cite{DrinfeldICM}) of a single bosonisation $\Lie{b} \rtimesdot \smun{\Lie{g}}{\mathrm{0}}$, which is in fact a subalgebra of the double-bosonisation.  Majid analyzed of the situation when $\rank \smun{\Lie{g}}{0}=\rank \Lie{g}-1$ for $\Lie{g}$ of finite type.  Sommerh\"{a}user has a similar construction (\cite{SommerhauserKacMoody}) and considered Kac--Moody algebras but only with $\smun{\Lie{g}}{0}=\Lie{h}$, the (generalized) Cartan subalgebra.

It is this construction that allows the inductive approach that was studied in detail in \cite{LieInduction} (and in \cite{QuantumLieInduction} for quantized enveloping algebras).  Our work in \cite{LieInduction} focused on finite-dimensional simple complex Lie algebras and showed that they have double-bosonisation decompositions related to deleting nodes from Dynkin diagrams.  We also considered necessary conditions on $\Lie{b}$ to obtain a simple Lie algebra after double-bosonisation.  In particular we obtained sufficiently many conditions to see that there could not be a finite-dimensional simple Lie algebra of type $E_{9}$, as an alternative method to the standard classification theorem.  In the course of the work in this paper, we have answered the natural follow-on question: which double-bosonisation \textit{does} give us the (infinite-dimensional) $E_{9}$?

In \cite{QuantumLieInduction}, we considered the corresponding situation for quantized enveloping algebras.  By means of Lusztig's root data approach (see for example \cite{LusztigBook}) we were able to prove corresponding results about double-bosonisation decompositions for all types, not just finite-type data.  Here we use the method of proof developed in \cite{QuantumLieInduction} to extend our earlier work for finite-dimensional Lie algebras to the infinite-dimensional cases, namely Kac--Moody algebras.

Our motivations here are two-fold.  First is the extension of our previous work, as described above, from the finite-dimensional Lie algebra situation to the infinite-dimensional Kac--Moody case.  Next is the identification of more examples of non-trivial braided-Lie bialgebras and double-bosonisation.  Previously-known examples of braided-Lie bialgebras have been mostly finite-dimensional; here we have shown that there is a large class of infinite-dimensional examples.  In particular, we have demonstrated a braided-Lie bialgebra structure on a well-known infinite-dimensional Lie algebra, the current algebra $k[u]\tensor \Lie{g}$ over a finite-dimensional simple Lie algebra $\Lie{g}$.  This braided coalgebra structure may be of interest in mathematical physics, where current algebras appear in quantum field theory.  Another direction for future work would be to tackle the other types of Kac--Moody algebra, in particular the twisted affine algebras.

The structure of this paper is as follows.  We begin by recalling in Section~\ref{s:prelims} the definitions of the structures we will need, namely those of a quasitriangular Lie bialgebra and a braided-Lie bialgebra.  We state the results of Majid (\cite{BraidedLie}) defining the double-bosonisation construction and its natural induced quasitriangular structure, when the input is quasitriangular.  We then recap the definition of Kac--Moody algebras associated to generalized Cartan matrices.  We will be interested in the untwisted affine Kac--Moody algebras and in particular their realizations as extensions of loop algebras over finite-dimensional simple Lie algebras.  

In Section~\ref{s:BLBAsassoctoKMAs}, we follow the pattern of \cite{QuantumLieInduction} by demonstrating that $\nat$-gradings of Lie bialgebras give rise to split projections.  Hence, via a theorem of Majid (\cite{BraidedLie}) corresponding to a result of Radford (\cite{Radford}) for Hopf algebras, one obtains a braided-Lie bialgebra associated to the grading as the kernel of the (split) projection.  Furthermore, one may reconstruct the original Lie bialgebra as a semidirect product Lie bialgebra of the kernel and image of the projection.

We then define a sub-root datum (Definition~\ref{subrootdatum}).  This is a way of formally encoding the idea of embedding Dynkin diagrams.  There is a natural $\nat$-grading of the negative Borel subalgebra of a Kac--Moody algebra associated to any choice of sub-root datum.  Hence by the general results we have a corresponding braided-Lie bialgebra (Theorem~\ref{bexists}).  Finally, it is clear from the definitions of the Kac--Moody algebras via a presentation and of this braided-Lie bialgebra that we can reconstruct the whole Kac--Moody algebra by a double-bosonisation (Proposition~\ref{gradedbraided}).  Equivalently, one may view Proposition~\ref{gradedbraided} as a generalisation of the usual triangular decomposition theorem.

In Section~\ref{s:BLBAsassoctoaffin} we focus on a particular case of the above results, namely the relationship between a finite-dimensional simple Lie algebra and the corresponding untwisted affine Kac--Moody algebra.  It is straightforward to see that if $\mathcal{C}$ is a root datum for an irreducible Cartan matrix of finite type then there is a natural root datum $\widetilde{\mathcal{C}}$ for the corresponding generalized Cartan matrix of untwisted affine type $\widetilde{C}$ and these fit together to give a sub-root datum $\mathcal{C} \subseteq_{\iota} \widetilde{\mathcal{C}}$ in the obvious way.  In terms of the Dynkin diagrams, we are simply considering the finite-type sub-diagram of the corresponding affine-type diagram, \ie the deletion of the extending node, usually labelled ``0''.  We call this sub-root datum the affinization sub-root datum associated to $\mathcal{C}$.

Hence by the general theory there is a $\nat$-grading of the Borel subalgebra of $L(\widetilde{C})$ and hence a braided-Lie bialgebra $\Lie{b}=\Lie{b}(\widetilde{\mathcal{C}},\mathcal{C},\iota)$.  By considering the loop algebra realization of $L(\widetilde{C})$ as $\hat{\Lie{L}}(C)$, we see that $\Lie{b}$ is isomorphic as a Lie algebra to a well-known Lie algebra, namely the current algebra over $L(C)$, $k[u]\tensor L(C)$.  However, we now see that the current algebra admits a non-trivial braided-Lie bialgebra structure.  We conclude with an example and some explicit formul\ae\ for the braided and unbraided cobrackets associated to $C=A_{2}$ (\ie $L(C)=\mathfrak{sl}_{3}(k)$, $L(\widetilde{C})\iso \hat{\Lie{L}}(A_{2})=\widetilde{\mathfrak{sl}}_{3}(k)$).

\subsection*{Acknowledgements}
The work in this paper has been completed during the author's research fellowship at Keble College, Oxford.  I would also like to thank the Mathematical Institute at Oxford for its provision of facilities.

I am grateful to the referee for several helpful comments and
corrections.

\section{Preliminaries}\label{s:prelims}

Throughout we work over a field $k=\bar{k}$, algebraically closed and of characteristic $0$.  We will use the following convention for the natural numbers: $\nat= \{ 0, 1, 2, \ldots \}$.  That is, for us $\nat$ is a monoid.

We will assume the reader is familiar with the theory of finite-dimensional semisimple Lie algebras and root systems, as can be found in \cite{Serre}, \cite{Humphreys} or \cite{CarterBook}, for example.  We shall recap the corresponding aspects of the theory for Kac--Moody algebras we need below.  Note that we use the term unqualified term ``Lie algebra'' inclusively, so that Kac--Moody algebras are considered to be examples of Lie algebras.

We use $\tau$ to denote the tensor product flip map, \eg \[ \tau:V\tensor W \to W\tensor V,\ \tau(v\tensor w)=w\tensor v\ \mbox{for all}\ v\in V,\ w\in W, \] on any appropriate pair of vector spaces.  The adjoint action of a Lie algebra $\Lie{g}$ on itself can be extended naturally to tensor products as follows.  For $x,y,z\in \Lie{g}$, \[ \ad{x}{y\tensor z}=\ad{x}{y}\tensor z +y\tensor \ad{x}{z}. \]  We use the term ad-invariant in the obvious way.  We adopt the Sweedler notation for elements of tensor products, namely we use upper or lower parenthesized indices to indicate the placement in the tensor product: $\sum a_{(1)} \tensor a_{(2)} \tensor a_{(3)} \in \Lie{g} \tensor \Lie{g} \tensor \Lie{g}$.  We will usually omit the summation sign.

The definition of a Lie bialgebra is originally due to Drinfel\cprime d (\cite{Drinfeld1}, \cite{DrinfeldICM}).  The idea is the same as that for Hopf algebras, where we have two structures dual to each other, compatible in a natural way.  It is worth commenting that Lie bialgebras form a richer class than Lie algebras: the choice of the cobracket, the dual structure to the bracket, is not usually unique.

\begin{definition}[{\cite{Drinfeld1}}]\label{Liebialg} A Lie bialgebra is $(\Lie{g},\Lbracket{\ }{\ },\delta)$ where 
\begin{enumerate} 
\item $(\Lie{g},\Lbracket{\ }{\ })$ is a Lie algebra, 
\item $(\Lie{g},\delta)$ is a Lie coalgebra, that is, $\delta:\Lie{g}\to \Lie{g}\tensor \Lie{g}$ satisfies 
	\[ \begin{array}{lll} \delta+\tau\circ \delta = 0 & & \mbox{\textup{(anticocommutativity)}} \\
			(\delta \tensor \id)\circ \delta +\mbox{\textup{cyclic}}=0 & & 
							\mbox{\textup{(co-Jacobi identity)}}
	\end{array} \]  (Here,\ \textup{``cyclic''}\ refers to cyclical rotations of the three tensor product factors in $\Lie{g}\tensor \Lie{g} \tensor \Lie{g}$.)
\item we have a cohomological compatibility condition: $\delta$ is a 1-cocycle in $Z_{\mathrm{ad}}^{1}(\Lie{g},\Lie{g}\tensor\Lie{g})$.  Explicitly, \[ \delta(\Lbracket{x}{y})=\ad{x}{\delta y}-\ad{y}{\delta x}. \] 
\end{enumerate}
\end{definition}

\noindent Examining this definition, we see that if $\Lie{g}$ is a finite-dimensional Lie bialgebra, then $(\dual{\Lie{g}},\dual{\delta},\dual{\Lbracket{\ }{\ }})$ is also a finite-dimensional Lie bialgebra.  Here, $\dual{\delta}$ and $\dual{\Lbracket{\ }{\ }}$ are the bracket and cobracket, respectively, given by dualisation.  Extending the strict notion of dualisation in the finite-dimensional case, we have the notion of dually paired Lie bialgebras.  We say Lie bialgebras $\Lie{f}$, $\Lie{g}$ are dually paired by a map $\ip{\ }{\ }:\Lie{f}\tensor \Lie{g} \to k$ if \[ \ip{\Lbracket{a}{b}}{g}=\ip{a\tensor b}{\delta g} \qquad \text{and} \qquad \ip{\delta a}{g\tensor h}=\ip{a}{\Lbracket{g}{h}} \] for $a,\, b\in \Lie{f}$, $g,\, h\in \Lie{g}$.  

In many cases, the cobracket $\delta$ arises as the coboundary of an element $r\in \Lie{g}\tensor\Lie{g}$.  (Explicitly, $\delta x=\ad{x}{r}$ for all $x\in \Lie{g}$.)  If $\delta$ further satisfies $(\id \tensor \delta)r=\Lbracket{\smun{r}{13}}{\smun{r}{12}}$, we say that $(\Lie{g},r)$ is a quasitriangular Lie bialgebra.  Here, we write $\smun{r}{12}=r^{(1)}\tensor r^{(2)} \tensor 1$, etc., with summation understood and the indices showing the placement in the triple tensor product $\Lie{g} \tensor \Lie{g} \tensor \Lie{g}$.  The bracket is taken in the common factor, so $\rule{0pt}{13pt}\Lbracket{\smun{r}{13}}{\smun{r}{12}}=\Lbracket{r^{(1)}}{{r'}^{(1)}} \tensor {r'}^{(2)} \tensor {r}^{(2)}$ with $r'$ a second copy of $r$. 

To construct a quasitriangular Lie bialgebra, it is sufficient to find an element $r\in \Lie{g} \tensor \Lie{g}$ satisfying the classical Yang--Baxter equation and with ad-invariant symmetric part.  Then we take the coboundary $\partial r$ for $\delta$.  The classical Yang--Baxter equation, in the Lie setting, is \[ \Lbbracket{r}{r} \defeq \Lbracket{\smun{r}{12}}{\smun{r}{13}} + \Lbracket{\smun{r}{12}}{\smun{r}{23}} + \Lbracket{\smun{r}{13}}{\smun{r}{23}} = 0. \]  The bracket $\Lbbracket{\ }{\ }$ is the Schouten bracket, the natural extension of the bracket to these tensor spaces.

We now consider the braided version of Lie bialgebras, as defined by Majid in \cite{BraidedLie}.  Here we consider the module category $\modcat{\Lie{g}}$ of a quasitriangular Lie bialgebra $\Lie{g}$ and objects in this category possessing a $\Lie{g}$-covariant Lie algebra structure.  We associate to these objects a braiding map generalising the usual flip, $\tau$.  If $\Lie{b}$ is a $\Lie{g}$-covariant Lie algebra in the category $\modcat{\Lie{g}}$, we define the infinitesimal braiding of $\Lie{b}$ to be the operator $\psi:\Lie{b}\tensor \Lie{b} \to \Lie{b}\tensor\Lie{b}$, $\psi(a\tensor b)=2r_{+} \laction (a\tensor b-b\tensor a)$ where $\laction$ is the left action of $\Lie{g}$ on $\Lie{b}$ extended to the tensor products.  In fact, $\psi$ is a $2$-cocycle in $Z_{\mathrm{ad}}^{2}(\Lie{b},\Lie{b}\tensor \Lie{b})$.

\begin{definition}[{\cite{BraidedLie}}]\label{blba} A braided-Lie bialgebra $(\Lie{b},\Lbracket{\ }{\ }_{\Lie{b}},\bdelta)$ is an object in $\modcat{\Lie{g}}$ satisfying the following conditions:  
\begin{enumerate} 
\item $(\Lie{b},\Lbracket{\ }{\ }_{\Lie{b}})$ is a $\Lie{g}$-covariant Lie algebra in the category, \ie $\Lbracket{\ }{\ }_{\Lie{b}}:\Lie{b} \tensor \Lie{b} \to \Lie{b}$ is a $\Lie{g}$-module map. 
\item $(\Lie{b},\bdelta)$ is a $\Lie{g}$-covariant Lie coalgebra in the category, \ie $\bdelta:\Lie{b} \to \Lie{b} \tensor \Lie{b}$ is a $\Lie{g}$-module map.
\item $\dform \bdelta=\psi$, or explicitly,
\[ \bdelta( [x,y] ) = \ad{x}{\bdelta y}-\ad{y}{\bdelta x}-\psi(x\tensor y). \] 
\end{enumerate}
\end{definition}

\noindent We may make a similar definition of dually paired braided-Lie bialgebras to that above for unbraided Lie bialgebras.  We can now state the theorem which provides the double-bosonisation construction of Majid.  Let $\Lie{g}$ be a quasitriangular Lie bialgebra.

\begin{theorem}[{\cite[Theorem 3.10]{BraidedLie}}]\label{dbos} For dually paired braided-Lie bialgebras $\Lie{b},\Lie{c}\in \modcat{\Lie{g}}$ the vector space $\Lie{b} \dsum \Lie{g} \dsum \Lie{c}$ has a unique Lie bialgebra structure $\dbos{b}{g}{\op{c}}$, the double-bosonisation, such that $\Lie{g}$ is a Lie sub-bialgebra, $\Lie{b},\Lie{\op{c}}$ are Lie subalgebras, and \[ \begin{array}{c} \Lbracket{\xi}{x}=\xi \laction x, \quad \Lbracket{\xi}{\phi}=\xi \laction \phi \\ \\  \Lbracket{x}{\phi} = x_{\underline{(1)}}\ip{\phi}{x_{\underline{(2)}}}+\phi_{\underline{(1)}}\ip{\phi_{\underline{(2)}}}{x}+2r_{+}^{(1)}\ip{\phi}{r_{+}^{(2)}\laction x} \\ \\ \delta x=\bdelta x+r^{(2)}\tensor r^{(1)}\laction x-r^{(1)}\laction x\tensor r^{(2)} \\ \\ \delta \phi=\bdelta \phi + r^{(2)}\laction \phi \tensor r^{(1)} - r^{(1)}\tensor r^{(2)}\laction \phi \end{array} \] \noindent for all $x\in \Lie{b}$, $\xi\in \Lie{g}$ and $\phi \in \Lie{c}$.  Here $\bdelta x=x_{\underline{(1)}}\tensor x_{\underline{(2)}}$. \end{theorem}

\noindent We note that Sommerh\"{a}user has a similar construction (\cite{SommerhauserKacMoody}).  Some different terminology is used there, notably `Yetter--Drinfel\cprime d Lie algebra' for `braided-Lie bialgebra'.  

The double-bosonisation is always quasitriangular when we take $\Lie{b}$ and $\Lie{c}=\dual{\Lie{b}}$ finite-dimensional, as we see from the following proposition.

\begin{proposition}[{\cite[Proposition 3.11]{BraidedLie}}]\label{doubleqtstr} Let $\Lie{b}\in \modcat{\Lie{g}}$ be a finite-dimensional braided-Lie bialgebra with dual $\dual{\Lie{b}}\!$.  Then the double-bosonisation $\dbos{b}{g}{\op{\dual{b}}}$ is quasitriangular with \[ r^{\mathrm{new}}=r+\sum_{a} f^a \tensor e_a \] where $\{ e_a \}$ is a basis of $\Lie{b}$ and $\{ f^a \}$ is a dual basis, and $r$ is the quasitriangular structure of $\Lie{g}$.  If $\Lie{g}$ is factorisable then so is the double-bosonisation.
\end{proposition}

\noindent We remark that for each finite-dimensional semisimple Lie algebra, there exists a canonical quasitriangular structure, associated to the Drinfel\cprime d--Sklyanin solution of the Yang--Baxter equation.

Next, we recall the definition of Kac--Moody algebras associated to generalized Cartan matrices, as introduced independently by Kac (\cite{Kac1967}) and Moody (\cite{Moody1967}) in 1967.  We mostly follow the presentation in \cite{CarterBook}, Chapters 14-19.  

A generalized Cartan matrix is an $n \cross n$ integer matrix $C$ such that $C_{ii}=2$ for all $i$, $C_{ij} \leq 0$ for $i\neq j$ and $C_{ij} = 0$ implies $C_{ji} =0$.  Let $I=\{ 1, \ldots , n \}$.  A triple $(H,\Pi,\ve{\Pi})$ is a minimal realization of $C$ over $k$ if $H$ is a $(2n-\rank C)$-dimensional $k$-vector space, $\Pi=\{ \alpha_{1} , \ldots , \alpha_{n} \} \subset \dual{H}$ and $\ve{\Pi} = \{ h_{1} , \ldots , h_{n} \} \subset H$ are linearly independent subsets in $\dual{H}$ and $H$ respectively and $\alpha_{j}(h_{i})=C_{ij}$ for all $i$.  A minimal realization allows us to define a Lie algebra $\widetilde{L}(C)$, by generators $\{ e_{1}, \ldots , e_{n}, f_{1} , \ldots , f_{n} , \tilde{h} \mid h\in H \}$ and relations
\begin{align*}
\tilde{x} & = \lambda \tilde{y} + \mu \tilde{z} & & \text{whenever}\ x=\lambda y + \mu z\ \text{holds in $H$}, \lambda, \mu \in k, \\
\Lbracket{\tilde{x}}{\tilde{y}} & = 0 & & \text{for all}\ x, y \in H, \\
\Lbracket{e_{i}}{f_{j}}& = \delta_{ij}\tilde{h_{i}} & & \text{for all}\ i,j\in I, \\
\Lbracket{\tilde{x}}{e_{i}} & = \alpha_{i}(x)e_{i} & & \text{for all}\ x\in H, i\in I, \\
\Lbracket{\tilde{x}}{f_{i}} & = -\alpha_{i}(x)f_{i} & & \text{for all}\ x\in H, i\in I.
\end{align*}
In fact, any two minimal realizations are isomorphic and hence $\widetilde{L}(C)$ is independent of the choice of minimal realization.  Let $\widetilde{H}$ denote the Abelian Lie subalgebra of $\widetilde{L}(C)$ generated by $\{ \tilde{h} \mid h\in H \}$, so $\widetilde{H} \iso H$ as $k$-vector spaces.  Then we define the Kac--Moody algebra $L(C)$ over $k$ associated to the generalized Cartan matrix $C$ to be the quotient of $\widetilde{L}(C)$ by the unique maximal ideal $J$ such that $J \intersection \widetilde{H} =(0)$.

A generalized Cartan matrix $C$ is called symmetrizable if there exists a non-singular diagonal matrix $D$ such that $CD$ is symmetric.  If $C$ is symmetrizable, then $L(C)$ has a presentation by the same generators as $\widetilde{L}(C)$ and relations as for $\widetilde{L}(C)$ together with the relations
\begin{eqnarray*}
(\mathrm{ad}\ e_{i})^{1-C_{ij}} e_{j} & = & 0 \\
(\mathrm{ad}\ f_{i})^{1-C_{ij}} f_{j} & = & 0
\end{eqnarray*}
where $\mathrm{ad}$ denotes the left adjoint action of $L(C)$ on itself via $\Lbracket{\ }{\ }$.  We note that if $C$ is in fact a Cartan matrix then $L(C)$ is the finite-dimensional semisimple Lie algebra with Cartan matrix $C$; the above presentation is the well-known Serre presentation (\cite{Serre}).

From now on, we will assume our generalized Cartan matrices are symmetrizable, unless otherwise stated.

We now recall the loop algebra construction and also the isomorphism of the output of this construction with the form given above by generators and relations.  This follows the exposition of Carter (\cite[Chapter~18]{CarterBook}).

Let $C$ be an irreducible Cartan matrix of finite type, of rank $l$, so that $L(C)$ is the associated finite-dimensional simple Lie algebra.  Then there exists a root $\theta=\sum_{i=1}^{l} a_{i}\alpha_{i}$ ($a_{i} \geq 0$) such that for any other root $\phi=\sum b_{i}\alpha_{i}$ we have $b_{i} \leq a_{i}$; we call $\theta$ the highest root of $L(C)$.  Then there is a coroot $h_{\theta}$ corresponding to $\theta$ and this may be written as $h_{\theta}=\sum_{i=1}^{l} c_{i}h_{i}$.  

Let $\widetilde{C}$ be the $(l+1) \cross (l+1)$ matrix with rows and columns indexed by $I=\{ 0, 1, \ldots, l \}$ such that
\begin{align*}
\widetilde{C}_{ij} & = C_{ij} & & \text{for}\ i, j \in \{ 1,\ldots ,l\}, \\
\widetilde{C}_{i0} & = -\sum_{j=1}^{l} a_{j}C_{ij} & & \text{for}\ i\in \{ 1, \ldots ,l \}, \\
\widetilde{C}_{0j} & = -\sum_{i=1}^{l} c_{i}C_{ij} & & \text{for}\ i\in \{ 1, \ldots, l \}, \\
\widetilde{C}_{00} & = 2. & &
\end{align*}
Then $\widetilde{C}$ is an affine Cartan matrix and if $X_{l}$ is the Dynkin type of $L(C)$, the type of $L(\widetilde{C})$ is $\widetilde{X}_{l}$.  The latter are the generalized Cartan matrices of untwisted affine type.

The loop algebra construction is as follows.  First, we take the tensor product of $L(C)$ with the $k$-vector space of Laurent polynomials in a variable $t$: $\Lie{L}(C) \defeq k[t, \inv{t}] \tensor L(C)$.  There is a unique Lie algebra structure on $\Lie{L}(C)$ given by $\Lbracket{p\tensor g}{q\tensor h}=pq\tensor \Lbracket{g}{h}$.  Next we make a $1$-dimensional central extension, setting $\tilde{\Lie{L}}(C) \defeq \Lie{L}(C) \dsum kc$.  The Lie algebra structure on $\tilde{\Lie{L}}(C)$ is \[ \Lbracket{(p\tensor g)+\lambda c}{(q\tensor h) + \mu c}=pq\tensor \Lbracket{g}{h}+\mathrm{Res}\left( \left( \frac{\dform p}{\dform t} \right)q \right)\ip{g}{h}c. \]  Here, $\mathrm{Res}$ is the residue function $\mathrm{Res}:k[t,\inv{t}]\to k$, $\mathrm{Res} \left( \sum \zeta_{i}t^{i} \right) = \zeta_{-1}$, and $\ip{\ }{\ }$ denotes the unique invariant bilinear form on $L(C)$ with $\ip{h_{\theta}}{h_{\theta}}=2$.  (We note that $\lambda$ and $\mu$ on the left in the above definition do not enter into the expression on the right, since $c$ is central in $\tilde{\Lie{L}}(C)$.)

Having made the central extension, we must also extend our algebra by a derivation.  This is done as follows.  Firstly, we see that $\Delta:\tilde{\Lie{L}}(C) \to \tilde{\Lie{L}}(C)$, $\Delta(p\tensor g+\lambda c)=t\frac{\dform p}{\dform t}\tensor g$ is a derivation, in the usual sense for linear maps.  Then we use $\Delta$ to define a Lie algebra structure on the vector space $\hat{\Lie{L}}(C)=\tilde{\Lie{L}}(C)\dsum kd$ by \[ \Lbracket{a+\lambda d}{b+\mu d}=\Lbracket{a}{b}+\lambda \Delta(b)-\mu \Delta(a) \] for $a,\, b\in \tilde{\Lie{L}}(C)$, $\lambda,\, \mu \in k$.  Putting everything together, this defines the following Lie bracket on $\hat{\Lie{L}}(C)$:
\begin{eqnarray*} \Lbracket{(p\tensor g)+\lambda c+\mu d}{(q\tensor h) + \sigma c + \tau d} & = & pq\tensor \Lbracket{g}{h}+\mathrm{Res}\left( \left( \frac{\dform p}{\dform t} \right)q \right)\ip{g}{h}c \\ & & \qquad + \mu t\frac{\dform q}{\dform t}\tensor h - \tau t\frac{\dform p}{\dform t}\tensor g. \end{eqnarray*}
This gives us a Lie algebra $\hat{\Lie{L}}(C)$, which we call the untwisted affine Kac--Moody algebra associated to the Cartan matrix $C$ of finite type.

Now we have two Lie algebras $L(\widetilde{C})$ and $\hat{\Lie{L}}(C)$ that we have claimed are the untwisted affine Kac--Moody algebras.  These algebras are in fact isomorphic and we now describe this isomorphism.  The method we sketch is not completely explicit: we use a result that says that if we can identify a generating set in a Kac--Moody algebra $\mathcal{A}$ that satisfies the relations defining $L(C)$ along with a realization of $C$, then in fact $\mathcal{A}\iso L(C)$.  The full details may be found in \cite[Theorem~18.5]{CarterBook}, from whence the following is taken.

To start, we consider the finite-dimensional simple Lie algebra $L(C)$ associated to an irreducible Cartan matrix $C$ of finite type and rank $l$.  Then $L(C)$ has a Serre presentation and in particular is generated by elements $E_{i}$, $F_{i}$ and $H_{i}$ for $1 \leq i \leq l$.  We let $\Lie{h}(C)$ denote the Cartan subalgebra of $L(C)$ with basis $\{ H_{i} \}$.  Therefore, in the derivation- and centrally-extended loop algebra $\hat{\Lie{L}}(C)$ we set
\begin{eqnarray*}
e_{i} & = & 1\tensor E_{i} \\
f_{i} & = & 1\tensor F_{i} \\
h_{i} & = & 1\tensor H_{i}
\end{eqnarray*}
for $i\in \{ 1,\ldots l \}$.  It is clear that these satisfy the required relations in $L(\widetilde{C})$ among themselves, since the submatrix of $\widetilde{C}$ given by deleting the first row and column is exactly $C$.

Next we need to define generators $e_{0}$, $f_{0}$, $h_{0}$.  Consider the root spaces $L(C)_{\theta}$, $L(C)_{-\theta}$, $\theta$ the highest root.  We choose $F_{0}\in L(C)_{\theta}$, $E_{0}\in L(C)_{-\theta}$ such that $\ip{F_{0}}{E_{0}}=1$ and $\omega(F_{0})=-E_{0}$ (where $\omega$ is the $L(C)$-automorphism such that $\omega(E_{i})=-F_{i}$, $\omega(F_{i})=-E_{i}$).

We set $e_{0}=t\tensor E_{0}$, $f_{0}=\inv{t}\tensor F_{0}$.  We need $h_{0}\in H=(1\tensor \Lie{h}(C))\dsum kc \dsum kd$ such that $\Lbracket{e_{0}}{f_{0}}=h_{0}$: we set $h_{0}=(1\tensor (-H_{\theta}))+c$.  Now we can also identify an element $\alpha_{0} \in \dual{H}$ such that $(H,\Pi,\ve{\Pi})$ with $\Pi=\{ h_{0} , h_{1}, \ldots , h_{l} \}$, $\ve{\Pi} = \{ \alpha_{0} , \alpha_{1} , \ldots , \alpha_{l} \}$ is a minimal realization of $\widetilde{C}$.  Putting all this together, we have that $\hat{\Lie{L}}(C)$ is isomorphic to $L(\widetilde{C})$ as a Lie algebra.  The isomorphism as Lie bialgebras is given by the same map: we simply push the cobracket on $L(\widetilde{C})$ over to $\hat{\Lie{L}}(C)$.

A Kac--Moody algebra $L(C)$ associated to a symmetrizable generalized Cartan matrix $C$ can be given the structure of a Lie bialgebra as follows.  We define $\delta:L(C) \to L(C) \tensor L(C)$ by $\delta e_{i} = \frac{d_{i}}{2}e_{i} \wedge \tilde{h}_{i}$ for all $i$, $\delta \tilde{x} = 0$ for all $x\in H$, $\delta f_{i} = \frac{d_{i}}{2} f_{i} \wedge \tilde{h}_{i}$ for all $i$, extended to the whole of $L(C)$ by the cocycle condition on $\delta$.  Here, the scalars $d_{i}$ are the non-zero entries of the diagonal matrix $D$ such that $CD$ is symmetric and $a\wedge b=a\tensor b-b\tensor a$ as usual.  We will only consider this Lie bialgebra structure on $L(C)$ and refer to $L(C)$ as a Kac--Moody bialgebra.

\section{Braided-Lie bialgebras associated to Kac--Moody bialgebras and double-bosonisation}\label{KacMoodyInduction}\label{s:BLBAsassoctoKMAs}

We will follow the pattern of our paper \cite{QuantumLieInduction}, first observing that braided-Lie bialgebras naturally arise whenever there is a $\nat$-grading of a Lie bialgebra.  We then introduce a datum which encodes embeddings of Lie bialgebras induced by inclusions of their Dynkin diagrams, see that there is an associated grading and hence a braided-Lie bialgebra.  We then see that we may reconstruct the full Lie bialgebra from the the subalgebra and the braided-Lie bialgebra by double-bosonisation.  We remark that by this method we reproduce earlier results of Majid (\cite[Proposition 4.5]{BraidedLie}) and ourselves (\cite[Proposition 3.3]{LieInduction}) as special cases.

An $\nat$-graded Lie bialgebra $\Lie{g}=\bigdsum_{n\in \nat} \Lie{g}_{n}$ is an $\nat$-graded Lie algebra in the usual sense whose Lie cobracket $\delta:\Lie{g} \to \Lie{g} \tensor \Lie{g}$ satisfies $\delta(\Lie{g}_{n}) \subseteq \bigdsum_{n=i+j} \Lie{g}_{i} \tensor \Lie{g}_{j}$.  Now any $\nat$-grading of a Lie bialgebra gives rise to a split Lie bialgebra projection, as described in the following lemma, which is easy to see.

\begin{lemma}\label{gradinggivessplproj} Let $\Lie{g}=\bigdsum_{n\in \nat} \Lie{g}_{n}$ be an $\nat$-graded Lie bialgebra.  Then $\Lie{g}_{0}$ is a Lie sub-bialgebra of $\Lie{g}$.  Let $\pi:\Lie{g} \onto \Lie{g}_{0}$ be defined by 
\[ \pi(\Lie{g}_{i}) = \begin{cases} \id|_{\Lie{g}_{0}} & \text{if}\ i=0 \\
                              0            & \text{otherwise.}
                \end{cases} \]
Then $\pi$ is a projection of $\nat$-graded Lie bialgebras, split by the inclusion $\iota:\Lie{g}_{0} \inj \Lie{g}$.  By this, we mean that $\pi$, $\iota$ are graded Lie bialgebra maps, such that $\pi$ is surjective, $\iota$ is injective and $\pi \circ \iota = \id_{\Lie{g}_{\scriptscriptstyle 0}}$ (the splitting condition).  $\Lie{g}_{0}$ is $\nat$-graded in the obvious way: $(\Lie{g}_{0})_{0}=\Lie{g}_{0}$, $(\Lie{g}_{0})_{i}=0$ $(i>0)$. \qed
\end{lemma}

Majid has given a version of Radford's theorem (\cite{Radford}), which associates to a split projection of Hopf algebras a braided Hopf algebra.  (These are also called Hopf algebras in braided categories and are the Hopf algebra version of a braided-Lie bialgebra, the former in fact being historically earlier.)  Majid's version for Lie bialgebras is as follows.

\begin{theorem}[{\cf \cite[Theorem 3.7]{BraidedLie}}]\label{LieRadfordMajid} Given any split projection of Lie bialgebras $\Lie{g} \genfrac{}{}{0pt}{}{\overset{\pi}{\onto}}{\underset{\iota}{\hookleftarrow}} \Lie{f}$, $\Lie{b} \defeq \Ker \pi$ is a braided-Lie bialgebra in the category of left $\Lie{f}$-modules with Lie algebra structure that of $\Ker \pi$ as a Lie subalgebra of $\Lie{g}$, $\Lie{f}$ acting by the left adjoint action and $\bdelta = (\id - \iota \circ \pi)^{\tensor 2} \circ \delta_{\Lie{g}}$.  Furthermore, $\Lie{b}$ admits a left $\Lie{f}$-coaction via $\beta:\Lie{b} \to \Lie{f} \tensor \Lie{b}$, $\beta=(\pi \tensor \id) \circ \delta_{\Lie{g}}$, and the action and coaction of $\Lie{f}$ are compatible so that $\Lie{b}$ is a Lie crossed module.  Finally, $\Lie{g} \iso \Lie{b} \rtimesdot \Lie{f}$, a simultaneous semidirect Lie algebra sum and coalgebra structure, called a (single) bosonisation.
\end{theorem}

We deduce that associated to any $\nat$-graded Lie bialgebra $\Lie{g}$, we have a graded braided-Lie bialgebra structure on $\Lie{b}=\Ker \pi=\bigdsum_{n > 0} \Lie{g}_{n}$.  The Lie algebra structure is just that of $\Lie{g}$; the novel feature is the braided cobracket.  We see from the expression for $\bdelta$ in the previous theorem that to calculate $\bdelta$ on any given element, we calculate $\delta$ and discard any terms where at least one tensor factor lies in degree $0$, \ie we keep only terms where both factors have degree at least one.  We remark that the Lie setting is considerably simpler to work in than that for Hopf algebras, since in the latter the resulting braided Hopf algebra is not the kernel of $\pi$ but is the space of coinvariants.  Here we obtain a grading on $\Lie{b}$ immediately, whereas this feature is harder to establish in the Hopf case.  The expression for the braided coproduct is also much harder to work with.

Next we define a sub-root datum, a formalisation of the idea of embedding Dynkin diagrams, adapting our definitions for the Hopf case.

\begin{definition} Let $C$ be a generalized Cartan matrix with columns indexed by a set $I$.  Let $(H,\Pi,\ve{\Pi})$ be a minimal realization of $C$, with $\Pi=\{ h_{i} \mid i\in I \}$, $\ve{\Pi}=\{ \alpha_{i} \mid i\in I \}$.  Then we say that $\mathcal{C}=(C,I,H,\Pi,\ve{\Pi})$ is a root datum associated to $C$.
\end{definition}

\begin{definition}\label{subrootdatum} Let $\mathcal{C}=(C,I,H,\Pi,\ve{\Pi})$, $\mathcal{C}'=(C',J,H',\Pi',\ve{(\Pi')})$ be two root data.  We say that $\mathcal{C}'$ is a sub-root datum of $\mathcal{C}$ via $\iota$ and write $\mathcal{C}' \subseteq_{\iota} \mathcal{C}$ if
\begin{enumerate}
\item $\iota:J \inj I$ is injective, 
\item $C'_{ij}=C_{\iota(i)\iota(j)}$ and
\item there exists an injective linear map $s:H' \inj H$ such that $s(h'_{i})=h_{\iota(i)}$ and $\dual{s}(\alpha'_{i})=\alpha_{\iota(i)}$, for all $h'_{i} \in \Pi'$, $\alpha'_{i}\in \ve{(\Pi')}$.  Here $\dual{s}:\dual{(H')} \inj \dual{H}$ is induced by $s$.
\end{enumerate}
So, we have $C'$ as a submatrix of $C$ (in a general sense), $H'$ identified with a subspace of $H$ and $\Pi'$, $\ve{(\Pi')}$ identified with subsets of $\Pi$, $\ve{\Pi}$ respectively.  It is clear that a sub-root datum induces an embedding of the Dynkin diagram associated to $C'$ in that associated to $C$.
\end{definition}

\noindent In particular, a sub-root datum induces an inclusion of Kac--Moody bialgebras:

\begin{lemma} Let $\mathcal{C'} \subseteq_{\iota} \mathcal{C}$ be a sub-root datum.  Then there is an injective Lie bialgebra homomorphism also denoted $\iota$ from $L(C')$ to $L(C)$.
\end{lemma}

\begin{proof} This follows from the presentations for $L(C')$ and $L(C)$, by defining $\iota$ on generators by $\iota(e_{i})=e_{\iota(i)}$, etc., and extending.  The form of the Lie cobracket makes it clear that this is a Lie bialgebra homomorphism.
\end{proof}

Now we can identify the $\nat$-graded Lie algebra which will yield the braided-Lie bialgebra we are seeking.  However, we cannot give $L(C)$ a suitable $\nat$-grading.  Instead, we want the analogue of the negative Borel subalgebra, defined as follows.

\begin{definition} Let $L(C)$ be the Kac--Moody algebra associated to a (not necessarily symmetrizable) generalized Cartan matrix $C$.  Define $B^{-}(C)$ to be the Lie subalgebra of $L(C)$ generated by the set $\{ f_{i},\, \tilde{h} \mid i\in I,\ h\in H \}$.  We will call $B^{-}(C)$ the negative Borel subalgebra of $L(C)$.
\end{definition}

Let $\mathcal{C'}=(C',J,H',\Pi',\ve{(\Pi')}) \subseteq_{\iota} \mathcal{C}$ be a sub-root datum of $C$ via $\iota$ and let $D=I\setminus \iota(J)$.  Denote by $\chi_{D}:I\to \{ 0, 1 \}$ the indicator function for $D$, \ie $\chi_{D}(i)=1$ if $i\in D$, $\chi_{D}(i)=0$ otherwise.

\begin{lemma} Associated to a sub-root datum $\mathcal{C'} \subseteq_{\iota} \mathcal{C}$, there is a $\nat$-grading of $B^{-}(C)$ defined by $\deg f_{i}=\chi_{D}(i)$ for all $i\in I$ and $\deg \tilde{h} = 0$.
\end{lemma}

\begin{proof} This follows by observing that the defining relations for $L(C)$ and hence $B^{-}(C)$ are homogeneous (\cf \cite[Section 1.5]{Kac}).
\end{proof}

\noindent Consequently, putting together Lemma~\ref{gradinggivessplproj} and Theorem~\ref{LieRadfordMajid} we obtain the following:

\begin{theorem}\label{bexists} Let $\mathcal{C'} \subseteq_{\iota} \mathcal{C}$ be a sub-root datum and let $B^{-}(C) \subset L(C)$ be the negative Borel subalgebra, $\nat$-graded by the grading coming from the sub-root datum.  Then there exists a braided-Lie bialgebra $\Lie{b}=\Lie{b}(\mathcal{C},\mathcal{C'},\iota)$ in the category of left $B^{-}(C)_{0}$-modules and left $B^{-}(C)_{0}$-comodules (in fact, Lie crossed modules).  We have $B^{-}(C) \iso \Lie{b} \rtimesdot B^{-}(C)_{0}$. \qed
\end{theorem}

It is easy to identify the structure of $B^{-}(C)_{0}$: it is generated by the set $\{ f_{j},\, \tilde{h} \mid j\in J,\ h\in H \}$ and is therefore, as a Lie algebra, a central extension of $B^{-}(C')$.  The Lie coalgebra structure is also clear.  Now although we have identified $\Lie{b}$ as a module for a central extension of $B^{-}(C')$, by the adjoint action, it is easy to verify that the $B^{-}(C')$ action extends to the adjoint action of $L(C')$.  Furthermore this is compatible with the central extension, so that $\Lie{b}$ is in fact a module for $L^{\sharp}(C')=\lgen e_{j},\, f_{j},\, \tilde{h} \mid j\in J,\ h\in H \rgen$.  

Also, we see that the $\nat$-grading on $B^{-}(C)$ induced by a sub-root datum may be modified to give a $\integ$-grading of the whole algebra $L(C)$ by defining $\deg e_{i}=\chi_{D}(i)$, $\deg f_{i}=-\chi_{D}(i)$, $\deg \tilde{x} = 0$.  Then $\Lie{b}$ may be identified as $\bigdsum_{n<0} L(C)_{n}$ (simply changing the sign of index for the grading) and the (graded) dual of $\Lie{b}$ can similarly be identified as $\bigdsum_{n > 0} L(C)_{n}$.  The degree $0$ part is exactly $L^{\sharp}(C')$.  Hence the following proposition is clear:

\begin{proposition}\label{gradedbraided} Let $\mathcal{C'} \subseteq_{\iota} \mathcal{C}$ be a sub-root datum and $L(C)$, $L(C')$ the associated Kac--Moody bialgebras.  Then we have the decomposition  \[ L(C) \iso \dbos{b}{\mathit{L^{\sharp}(C')}}{\op{\dual{b}}}. \] 
\end{proposition}

\begin{proof} This follows from our discussions above and from examining the presentation for $L(C)$.  An alternative more formal argument follows the line taken in \cite{QuantumLieInduction}, by constructing $L(C)$ in two ways as quotients of the Drinfel\cprime d doubles of $B^{-}(C)$ and $\Lie{b} \rtimesdot B^{-}(C)_{0}$.  (It is a general feature that double-bosonisation has a realization as a quotient of the double of a (single) bosonisation.)
\end{proof}

\begin{jegrmk}
\textit{Recalling the definition of double-bosonisation, we see that we ought to have $L^{\sharp}(C')$ quasitriangular but this will not be true if this algebra is infinite-dimensional.  However, this is a technicality: $L^{\sharp}(C')$ is pseudo-quasitriangular---the problem is one of taking infinite sums---so that the necessary axioms hold in a formal sense and in fact properties such as the local finiteness of the adjoint action mean that only ``finitely much'' of the sum is needed at any one time.  Indeed, the fact that we already know that $L(C)$ admits a genuine Lie bialgebra structure tells us that the double-bosonisation is sufficiently well-defined and that is why we have described this as a ``decomposition''.}

\textit{For a better resolution to this issue, one should use an appropriate version of Majid's weakening of the axioms of a quasitriangular structure for Hopf algebras, where the element $\mathcal{R}$ is replaced by certain maps.  In \cite{QuantumLieInduction}, we reformulated of Majid's definition to that of a so-called weak quasitriangular system and the key point is that double-bosonisation is still well-defined with respect to such a system for the middle algebra.  The corresponding Lie version would then give the right formalism to state the above proposition more precisely.}
\end{jegrmk}

This proposition extends the previous results of Sommerh\"{a}user (\cite{SommerhauserKacMoody}) (for $C'=C$), Majid (\cite{BraidedLie}) (the case of $C$ an $l \cross l$ irreducible Cartan matrix of finite type, $C'$ $(l-1) \cross (l-1)$ and irreducible) and our own earlier result (\cite{LieInduction}) (not assuming $C'$ irreducible), although not quite achieving the generality of the quantum case (\cite{QuantumLieInduction}).

\section{Braided-Lie bialgebras associated to affinization}\label{s:BLBAsassoctoaffin}

In this section, we focus on untwisted affine Kac--Moody bialgebras.  We have seen that these Kac--Moody bialgebras have realizations as extensions of loop algebras of finite-dimensional simple Lie (bi-)algebras.  We want to show that if $C$ is an irreducible Cartan matrix of finite type and $\widetilde{C}$ is the affine generalized Cartan matrix associated to $C$ then there is a sub-root datum $\mathcal{C} \subseteq_{\iota} \widetilde{\mathcal{C}}$ and hence an associated braided-Lie bialgebra.  Our goal is to analyse this braided-Lie bialgebra structure.

We first fix $\mathcal{C}=(C,I,H,\Pi,\ve{\Pi})$ a root datum for $C$ with $I=\{ 1,\ldots , l \}$, $\Pi=\{ h_{i} \mid i \in I \}$, $\ve{\Pi}=\{ \alpha_{i} \mid i\in I \}$.  Then we have $\widetilde{\mathcal{C}}=(\widetilde{C},\widetilde{I},\widetilde{H},\widetilde{\Pi},\ve{\widetilde{\Pi}})$ with $\widetilde{I}=\{ 0, 1, \ldots , l \}$, $\widetilde{\Pi}=\{ h_{0} \} \union \Pi$, $\ve{\widetilde{\Pi}}=\{ \alpha_{0} \} \union \ve{\Pi}$.  The sub-root datum $\mathcal{C} \subseteq_{\iota} \widetilde{\mathcal{C}}$ is entirely natural: we simply take $\iota$ to be the map $\iota:\{ 1,\ldots , l\} \to \{0, 1, \ldots , l \}$, $\iota(i)=i$ and then the maps $s$ and $\dual{s}$ in the definition of a sub-root datum are $h_{i} \mapsto h_{i}$ and $\alpha_{i} \mapsto \alpha_{i}$, for $1\leq i \leq l$.  We will call $\mathcal{C} \subseteq_{\iota} \widetilde{\mathcal{C}}$ the affinization sub-root datum associated to $\mathcal{C}$.

Theorem~\ref{bexists} tells us that we have a braided-Lie bialgebra $\Lie{b}=\Lie{b}(\widetilde{\mathcal{C}},\mathcal{C},\iota)$ in the category of left $B^{-}(C)$-modules.  However we want to know the structure of $\Lie{b}$ more explicitly as the kernel of a Lie bialgebra projection associated to a grading on the loop realization of $L(\widetilde{C})$.  

We note that there is a natural $\integ$-grading on the loop algebra $\Lie{L}(C)$: the obvious one, namely $\deg (t^{i}\tensor g)=i$ for $g\in L(C)$.  By setting $\deg c = \deg d = 0$ and examining the definition of the bracket on $\hat{\Lie{L}}(C)$, we see that this extends to a $\integ$-grading on the whole of $\hat{\Lie{L}}(C)$.  Now recall the isomorphism of $L(\widetilde{C})$ with $\hat{\Lie{L}}(C)$ and observe that the generators $e_{i}$, $h_{i}$ and $f_{i}$ for $1\leq i \leq l$ are mapped to elements in $\hat{\Lie{L}}(C)$ of degree $0$, as is $h_{0}$.  Furthermore, $e_{0}$ and $f_{0}$ are mapped to elements of degree $1$ and $-1$, respectively.  The projection $\pi$ that gives us the braided-Lie bialgebra $\Lie{b}$ we seek is the bialgebra projection $\pi: \bigdsum_{n\leq 0} \hat{\Lie{L}}(C)_{n} \onto \hat{\Lie{L}}(C)_{0}=B^{-}(C)\dsum kc \dsum kd$ which is the identity on homogeneous elements of degree $0$ and is zero on homogeneous elements of non-zero degree.  This allows us to identify $\Lie{b}$ explicitly.

\begin{corollary}
Let $\Lie{b}=\Lie{b}(\widetilde{\mathcal{C}},\mathcal{C},\iota)$ be the braided-Lie bialgebra associated to the affinization sub-root datum $\mathcal{C} \subseteq_{\iota} \widetilde{\mathcal{C}}$.  Then $\Lie{b}=\bigdsum_{n<0} \hat{\Lie{L}}(C)\iso \inv{t}k[\inv{t}]\tensor L(C)\iso k[u]\tensor L(C)$ as a Lie algebra. \qed
\end{corollary}

We remark that the Lie algebra $k[u]\tensor L(C)$ is well-known as the current algebra associated to $L(C)$.  The novel aspect here is the braided cobracket on the current algebra, making it a braided-Lie bialgebra.

To end, we give an example and some explicit formul\ae\ for the braided cobrackets arising in this way.

\begin{example*} We take $C={A}_{2}$ so $L(C)=\mathfrak{sl}_{3}(k)$ and $\widetilde{C}=\widetilde{A}_{2}$.  We have used the isomorphism of the two forms $L(\widetilde{A}_{2})$ and $\hat{\Lie{L}}(A_{2})$ to calculate the cobracket on elements of a basis for $\hat{\Lie{L}}(A_{2})$.  The basis we use is $\{ E_{1}, E_{2}, E_{12}, H_{1}, H_{2}, F_{1}, F_{2}, F_{21} \}$ where $E_{1}$, $E_{2}$, $H_{1}$, $H_{2}$, $F_{1}$ and $F_{2}$ are generators and $E_{12}=\Lbracket{E_{1}}{E_{2}}$ and $F_{21}=\Lbracket{F_{2}}{F_{1}}$.  Since $\widetilde{C}$ is symmetric, we have $d_{i}=1$ for all $i$.  Then one extends the definition of $\delta$ on generators to the whole of $\hat{\Lie{L}}(A_{2})$ via the cocycle condition (Definition~\ref{Liebialg} \textit{iii})).  We use this information to calculate the braided cobracket on $\dual{\Lie{b}}=\dual{(\Lie{b}(\widetilde{\mathcal{A}}_{2},\mathcal{A}_{2},\iota))}$, where $\mathcal{A}_{2}$ is the root datum corresponding to the usual root system for $\mathfrak{sl}_{3}(k)$ and $\widetilde{\mathcal{A}}_{2}$ similarly for $\widetilde{\mathfrak{sl}}_{3}$.  We work with the dual of $\Lie{b}$ simply for notational convenience, because $\Lie{b}$ itself consists of negative degree elements whereas $\dual{\Lie{b}}$ lives in positive degrees.  It is of course trivial to convert from $\dual{\Lie{b}}$ back to $\Lie{b}$.

We now give the cobracket on basis elements of $\hat{\Lie{L}}(A_{2})$ of degree at least $1$ (\ie $i>0$).
{\allowdisplaybreaks
\begin{eqnarray*} 
\delta(t^{i} \tensor E_{1}) & = & \frac{1}{2} (t^{i} \tensor E_{1}) \wedge (ic+1\tensor H_{1}) + \displaystyle\sum_{j=0}^{i-1} (t^j \tensor E_{1}) \wedge (t^{i-j} \tensor H_{1}) \\ & & \qquad  - \displaystyle\sum_{j=0}^{i-1} (t^j \tensor E_{12}) \wedge (t^{i-j} \tensor F_{2}) \\
\delta(t^{i} \tensor E_{2}) & = & \frac{1}{2}(t^{i} \tensor E_{2}) \wedge (ic+1\tensor H_{2}) + \displaystyle\sum_{j=0}^{i-1} (t^j \tensor E_{2}) \wedge (t^{i-j} \tensor H_{2}) \\ & & \qquad  + \displaystyle\sum_{j=0}^{i-1} (t^j \tensor E_{12}) \wedge (t^{i-j} \tensor F_{1}) \\
\delta(t^{i} \tensor E_{12}) & = & \frac{1}{2}(t^{i} \tensor E_{12}) \wedge (ic+1\tensor (H_{1}+H_{2})) + \displaystyle\sum_{j=0}^{i-1} (t^j \tensor E_{12}) \wedge (t^{i-j} \tensor (H_{1}+H_{2})) \\ & & \qquad  + \displaystyle\sum_{j=0}^{i} (t^j \tensor E_{2}) \wedge (t^{i-j} \tensor E_{1}) \\
\delta(t^i \tensor H_{1}) & = & \frac{1}{2}(t^i \tensor H_{1}) \wedge (ic) - 2 \displaystyle\sum_{j=0}^{i-1} (t^j \tensor E_{1}) \wedge (t^{i-j} \tensor F_{1}) \\ & & \qquad + \displaystyle\sum_{j=0}^{i-1} (t^j \tensor E_{2}) \wedge (t^{i-j} \tensor F_{2}) - \displaystyle\sum_{j=0}^{i-1} (t^j \tensor E_{12}) \wedge (t^{i-j} \tensor F_{21}) \\
\delta(t^i \tensor H_{2}) & = & \frac{1}{2}(t^i \tensor H_{2}) \wedge (ic) - 2 \displaystyle\sum_{j=0}^{i-1} (t^j \tensor E_{2}) \wedge (t^{i-j} \tensor F_{2}) \\ & & \qquad + \displaystyle\sum_{j=0}^{i-1} (t^j \tensor E_{1}) \wedge (t^{i-j} \tensor F_{1}) - \displaystyle\sum_{j=0}^{i-1} (t^j \tensor E_{12}) \wedge (t^{i-j} \tensor F_{21}) \\
\delta(t^i \tensor F_{1}) & = & \frac{1}{2}(t^i \tensor F_{1}) \wedge (ic+1\tensor H_{1}) - \displaystyle\sum_{j=1}^{i} (t^{j} \tensor F_{1}) \wedge (t^{i-j} \tensor H_{1}) \\ & & \qquad + \displaystyle\sum_{j=1}^{i} (t^j \tensor F_{21}) \wedge (t^{i-j} \tensor E_{2}) \\
\delta(t^i \tensor F_{2}) & = & \frac{1}{2}(t^i \tensor F_{2}) \wedge (ic+1\tensor H_{2}) - \displaystyle\sum_{j=1}^{i} (t^{j} \tensor F_{2}) \wedge (t^{i-j} \tensor H_{2}) \\ & & \qquad - \displaystyle\sum_{j=1}^{i} (t^j \tensor F_{21}) \wedge (t^{i-j} \tensor E_{1}) \\
\delta(t^{i} \tensor F_{21}) & = & \frac{1}{2}(t^{i} \tensor F_{21}) \wedge (ic+1\tensor (H_{1}+H_{2})) - \displaystyle\sum_{j=1}^{i} (t^j \tensor F_{21}) \wedge (t^{i-j} \tensor (H_{1}+H_{2})) \\ & & \qquad  + \displaystyle\sum_{j=1}^{i-1} (t^j \tensor F_{1}) \wedge (t^{i-j} \tensor F_{2}) 
\end{eqnarray*}
}

Via the theorems of Section~\ref{KacMoodyInduction}, we obtain the braided cobracket $\bdelta$ from $\bdelta=(\id - \pi)^{\tensor 2} \circ \delta$, where $\pi$ is the projection that kills all basis elements of non-zero degree.  Hence $\id - \pi$ kills basis elements of degree $0$ and $(\id - \pi)^{\tensor 2}$ keeps only those terms in $\delta$ with both tensor factors having non-zero degree.  This gives us the braided cobracket below, on $\dual{\Lie{b}}=\bigdsum_{n>0} \hat{\Lie{L}}(A_{2})_{n} \iso k[u]\tensor \mathfrak{sl}_{3}(k)$.

For $i>0$,
{\allowdisplaybreaks
\begin{eqnarray*}
\bdelta(t^{i} \tensor E_{1}) & = & \displaystyle\sum_{j=1}^{i-1} (t^j \tensor E_{1}) \wedge (t^{i-j} \tensor H_{1}) - \displaystyle\sum_{j=1}^{i-1} (t^j \tensor E_{12}) \wedge (t^{i-j} \tensor F_{2}) \\
\bdelta(t^{i} \tensor E_{2}) & = & \displaystyle\sum_{j=1}^{i-1} (t^j \tensor E_{2}) \wedge (t^{i-j} \tensor H_{2}) + \displaystyle\sum_{j=1}^{i-1} (t^j \tensor E_{12}) \wedge (t^{i-j} \tensor F_{1}) \\
\bdelta(t^{i} \tensor E_{12}) & = & \displaystyle\sum_{j=1}^{i-1} (t^j \tensor E_{12}) \wedge (t^{i-j} \tensor (H_{1}+H_{2})) + \displaystyle\sum_{j=1}^{i-1} (t^j \tensor E_{2}) \wedge (t^{i-j} \tensor E_{1}) \\
\bdelta(t^i \tensor H_{1}) & = & \displaystyle\sum_{j=1}^{i-1} (t^j \tensor E_{2}) \wedge (t^{i-j} \tensor F_{2}) -2 \displaystyle\sum_{j=1}^{i-1} (t^j \tensor E_{1}) \wedge (t^{i-j} \tensor F_{1}) \\ & & \qquad - \displaystyle\sum_{j=1}^{i-1} (t^j \tensor E_{12}) \wedge (t^{i-j} \tensor F_{21}) \\
\bdelta(t^i \tensor H_{2}) & = & \displaystyle\sum_{j=1}^{i-1} (t^j \tensor E_{1}) \wedge (t^{i-j} \tensor F_{1}) - 2 \displaystyle\sum_{j=1}^{i-1} (t^j \tensor E_{2}) \wedge (t^{i-j} \tensor F_{2}) \\ & & \qquad - \displaystyle\sum_{j=1}^{i-1} (t^j \tensor E_{12}) \wedge (t^{i-j} \tensor F_{21}) \\
\bdelta(t^i \tensor F_{1}) & = & - \displaystyle\sum_{j=1}^{i-1} (t^{j} \tensor F_{1}) \wedge (t^{i-j} \tensor H_{1}) + \displaystyle\sum_{j=1}^{i-1} (t^j \tensor F_{21}) \wedge (t^{i-j} \tensor E_{2}) \\
\bdelta(t^i \tensor F_{2}) & = &  - \displaystyle\sum_{j=1}^{i-1} (t^{j} \tensor F_{2}) \wedge (t^{i-j} \tensor H_{2}) - \displaystyle\sum_{j=1}^{i-1} (t^j \tensor F_{21}) \wedge (t^{i-j} \tensor E_{1}) \\
\bdelta(t^{i} \tensor F_{21}) & = & - \displaystyle\sum_{j=1}^{i-1} (t^j \tensor F_{21}) \wedge (t^{i-j} \tensor (H_{1}+H_{2})) + \displaystyle\sum_{j=1}^{i-1} (t^j \tensor F_{1}) \wedge (t^{i-j} \tensor F_{2})
\end{eqnarray*}
}
\end{example*}

\small

\bibliographystyle{halpha}
\bibliography{references}\label{references}
\addcontentsline{toc}{chapter}{References}

\normalsize

\end{document}